%% file: stability.tex
\documentclass[11pt,twoside]{article}
\usepackage{amsmath,amsfonts,amscd,amsthm}
\newtheorem{thm}{Theorem}[section]
\newtheorem{prop}[thm]{Proposition}

\newtheorem{cor}[thm]{Corollary}

\theoremstyle{definition}
\newtheorem{defin}[thm]{Definition}

\theoremstyle{remark}
\newtheorem{rem}[thm]{Remark}
\numberwithin{equation}{section}

\input macros   %abreviations diverses
\begin{document}

\title{Stability of the vanishing of the $\opa_b$-cohomology under small
horizontal perturbations of the CR structure in compact abstract
$q$-concave CR manifolds}

\author{Christine LAURENT-THI\'{E}BAUT}

\date{Pr\'{e}publication de l'Institut Fourier n$^\circ$~xxx (2000)\\
\vspace{5pt} {\tt
http://www-fourier.ujf-grenoble.fr/prepublications.html }}
\date{}
\maketitle

\ufootnote{\hskip-0.6cm  {\it A.M.S. Classification}: 32V20.
\newline {\it Key words}: CR structure, Homotopy formulas, Tangential
Cauchy Riemann equation, Vanishing theorems.}

\bibliographystyle{amsplain}
%\setcounter{page}{-1}
%\input resume

%\tableofcontents
\input anis0
\input anis1
\input anis3
\input anis4

\bibliography{biblio}
\vespa\vespa
\begin{flushright}
  \begin{minipage}[t]{8cm}

Universit\'e de Grenoble

   Institut Fourier

UMR 5582  CNRS/UJF

   BP 74

38402 St Martin d'H\`eres Cedex

France

Christine.Laurent@ujf-grenoble.fr

  \end{minipage}
\end{flushright}

\enddocument

\end

%% file: macros.tex
\providecommand\ufootnote[1]{{\let\thefootnote\relax\footnote[0]{#1}}}

\newcommand{\vespa}{\vspace{1em}}

\newcommand{\ac}{\mathcal A}
\newcommand{\bc}{\mathcal B}

\newcommand{\cc}{\mathcal C}

\newcommand{\ec}{\mathcal E}

\newcommand{\cs}{\mathcal S}

\newcommand{\rb}{\mathbb R}
\newcommand{\nb}{\mathbb N}
\newcommand{\mb}{\mathbb M}
\newcommand{\cb}{\mathbb C}

\newcommand{\ci}{{\mathcal C}^\infty}

\newcommand{\Cal}{\mathcal}

\newcommand{\ol}{\overline}

\newcommand{\pa}{\partial}
\newcommand{\opa}{\ol\partial}
\newcommand{\wt}{\widetilde}

 \DeclareMathOperator{\im}{Im}
\DeclareMathOperator{\ke}{Ker}

%% file: anis0.tex
The tangential Cauchy-Riemann equation is one of the main tools 
in CR analysis 
and its properties are deeply related to the geometry of CR manifolds, in 
particular the complex tangential 
directions are playing an important role. For example it was noticed
by Folland and Stein \cite{FoSt}, when they studied the tangential
Cauchy-Riemann operator on the Heisenberg group and more generally
on strictly pseudoconvex real hypersurfaces of $\cb^n$, that one get
better estimates in the complex tangential directions. Therefore in the study of the stability properties for the  tangential Cauchy-Riemann equation under perturbations of the CR structure it
seems natural to consider perturbations which
preserve the complex tangent vector bundle. Such perturbations can
be represented as graphs in the complex tangent vector bundle over
the original CR structure, they are defined by $(0,1)$-forms with
values in the holomorphic tangent bundle. We call them horizontal
perturbations.

We consider compact abstract CR manifolds and integrable
perturbations of their CR structure which preserve their complex
tangent vector bundle. Since the Levi form of a CR manifold depends
only on its complex tangent vector bundle, such perturbations will
preserve the Levi form and hence the concavity properties of the
manifold which are closely related to the $\opa_b$-cohomology. For
example it is well known that, for $q$-concave compact CR manifolds
of real dimension $2n-k$ and CR-dimension $n-k$, the
$\opa_b$-cohomology groups are finite dimensional in bidegree
$(p,r)$, when $1\le p\le n$ and $1\le r\le q-1$ or $n-k-q+1\le r\le
n-k$. Therefore, for a $q$-concave compact CR manifold, this
finiteness property is stable by horizontal perturbations of the CR
structure.

In this paper we are interested in the stability of the vanishing of
the $\opa_b$-cohomology groups after horizontal perturbations of the
CR structure.

Let $\mb=(\mb,H_{0,1}\mb)$ be an abstract compact CR manifold of
class $\ci$, of real dimension $2n-k$ and CR dimension $n-k$, and
$\widehat\mb=(\mb,\widehat H_{0,1}\mb)$ another abstract compact CR
manifold such that $\widehat H_{0,1}\mb$ is a smooth integrable
horizontal perturbation of $H_{0,1}\mb$, then $\widehat H_{0,1}\mb$
is defined by a smooth form $\Phi\in\ci_{0,1}(\mb,H_{1,0}\mb)$. We
denote by $\opa_b$ the tangential Cauchy-Riemann operator associated
to the CR structure $H_{0,1}\mb$ and by $\opa_b^\Phi$ the tangential
Cauchy-Riemann operator associated to the CR structure $\widehat
H_{0,1}\mb$.

The smooth $\opa_b$-cohomology groups on $\mb$ in bidegree $(0,r)$
and $(n,r)$ are defined for $1\le r\le n-k$ by:
$$H^{0,r}(\mb)=\{f\in\ci_{0,r}(\mb)~|~\opa_b f=0\}/\opa_b (\ci_{0,r-1}(\mb))$$
and
$$H^{n,r}(\mb)=\{f\in\ci_{n,r}(\mb)~|~\opa_b f=0\}/ \opa_b (\ci_{n,r-1}(\mb)).$$

If $f$ is a smooth differential form of degree $r$, $1\le r\le n-k$,
on $\mb$, we denote by $f_{r,0}$ its projection on the space
$\ci_{r,0}(\mb)$ of $(r,0)$-forms for the CR structure $H_{0,1}\mb$. Note that if $r\geq n+1$ then $f_{r,0}=0$.

The smooth $\opa_b^\Phi$-cohomology groups on $\widehat \mb$ in
bidegree $(0,r)$ and $(n,r)$ are defined for $1\le r\le n-k$ by:
$$H^{0,r}_\Phi(\widehat\mb)=\{f\in\ci_{0,r}(\widehat\mb)~|~f_{r,0}=0,~\opa_b^\Phi f=0\}/\opa_b^\Phi (\ci_{o,r-1}(\widehat\mb))$$
and
$$H^{n,r}_\Phi(\widehat\mb)=\{f\in\ci_{n,r}(\widehat\mb)~|~\opa_b^\Phi f=0\}/ \opa_b^\Phi (\ci_{n,r-1}(\widehat\mb)).$$

In this paper the following stability result is proved:
\begin{thm}\label{stabilite}
Assume $\mb$ is $q$-concave, there exists then a sequence
$(\delta_l)_{l\in\nb}$ of positive real numbers such that, if
$\|\Phi\|_{l}<\delta_l$ for each $l\in\nb$,

(i)  $H^{p,r-p}(\mb)=0$ for all $1\leq p\leq r$ implies
$H^{0,r}_\Phi(\widehat\mb)=0$ when $1\leq r\leq q-2$, in the
abstract case, and when $1\leq r\leq q-1$, if $\mb$ is locally
embeddable,

(ii) $H^{n-p,r+p}(\mb)=0$ for all $0\le p\le n-k-r$ implies
$H^{n,r}_\Phi(\widehat\mb)=0$ when $n-k-q+1\le r\le n-k$.
\end{thm}

We also prove the stability of the solvability of the tangential
Cauchy-Riemann equation with sharp anisotropic regularity (cf.
Theorem \ref{annul}).

Note that, when both CR manifolds $\mb$ and $\widehat \mb$ are
embeddable in the same complex manifold (i.e. in the embedded case),
a $(0,r)$-form $f$ for the new CR structure $\widehat H_{0,1}\mb$ is
also a $(0,r)$-form for the original CR structure $H_{0,1}\mb$ and
hence the condition $f_{r,0}=0$ in the definition of the cohomology
groups $H^{0,r}_\Phi(\widehat\mb)$ is automatically fulfilled. In that case,
Polyakov proved in \cite{Poly1} global homotopy formulas for a
family of CR manifolds in small degrees which immediately imply the
stability of the vanishing of the $\opa_b$-cohomology groups of
small degrees. In his paper he does not need the perturbation to
preserve the complex tangent vector bundle, but his estimates are
far to be sharp.

Moreover Polyakov \cite{Poly2} proved also that if a generically
embedded compact CR manifold $M\subset X$ is at least $3$-concave
and satisfies $H^{0,1}(M,T'X_{|_M})=0$, then small perturbations of
the CR structure are still embeddable in the same manifold $X$. From
this result and the global homotopy formula from \cite{Poly1} one can derive some stability results for the vanishing
of $\opa_b$-cohomology groups of small degrees without hypothesis on the embeddability
a priori of the perturbed CR structure.

The main interest of our paper is that we do not assume the CR
manifolds $\mb$ and $\ol\mb$ to be embeddable. In the case where the original CR
manifold $\mb$ is embeddable it covers the case where it is unknown if the
perturbed CR structure is embeddable in the same manifold as the
original one, for example when the manifold $\mb$ is only $2$-concave.
 Finally we also reach
the case of the $\opa_b$-cohomology groups of large degrees, which,
even in the embedded case, cannot be deduced from the works of
Polyakov.

The main tool in the proof of the stability of the vanishing of the
$\opa_b$-cohomology groups is a fixed point theorem which is derived
from global homotopy formulas with sharp anisotropic estimates. Such
formulas are proved in \cite{ShWa}, in the abstract case, by using
the $L^2$ theory for the $\Box_b$ operator and in \cite{Laanis}, in
the locally embeddable case, by first proving that the integral
operators associated to the kernels built in \cite{BaLa} satisfies
sharp anisotropic estimates, which implies local homotopy formulas
with sharp anisotropic estimates, and then by using the
globalization method from \cite{LaLeglobal}and \cite{BroLaLe}.

%% file: anis1.tex
\section{CR Structures}\label{s1}

Let $\mb$ be a $\cc^l$-smooth, $l\geq 1$, paracompact differential manifold, we
denote by $T\mb$ the tangent bundle of $\mb$ and by $T_\cb
\mb=\cb\otimes T\mb$ the complexified tangent bundle.

\begin{defin}\label{structure}
An \emph {almost CR structure} on $\mb$ is a subbundle $H_{0,1}\mb$ of
$T_\cb \mb$ such that $H_{0,1}\mb\cap \ol{H_{0,1}\mb}=\{0\}$.

If the almost CR structure is integrable, i.e. for all
$Z,W\in\Gamma(\mb,H_{0,1}\mb)$ then $[Z,W]\in\Gamma(\mb,H_{0,1}\mb)$, then it 
is called a \emph {CR structure}.

If $H_{0,1}\mb$ is a CR structure, the pair $(\mb,H_{0,1}\mb)$ is
called an \emph {abstract CR manifold}.

The \emph{CR dimension} of $\mb$ is defined by CR-dim
$\mb=\rm{rk}_\cb~H_{0,1}\mb$.
\end{defin}

We set $H_{1,0}\mb=\ol{H_{0,1}\mb}$ and we denote by $H^{0,1}\mb$
the dual bundle $(H_{0,1}\mb)^*$ of $H_{0,1}\mb$.

Let $\Lambda^{0,q}\mb=\bigwedge^q(H^{0,1}\mb)$, then
$\cc_{0,q}^s(\mb)=\Gamma^s(\mb,\Lambda^{0,q}\mb)$ is called the
space of $(0,q)$-forms of class $\cc^s$, $0\leq s\leq l$, on $\mb$.

We define $\Lambda^{p,0}\mb$ as the space of forms of degree $p$ that 
annihilate any $p$-vector on $\mb$ that has more than one factor contained 
in $H_{0,1}\mb$. Then $\cc_{p,q}^s(\mb)=\cc_{0,q}^s(\mb,\Lambda^{p,0}\mb)$ 
is the space of $(0,q)$-forms of class $\cc^s$ with values in 
$\Lambda^{p,0}\mb$.  

If the almost CR structure is a CR structure, i.e. if it is
integrable, and if $s\geq 1$, then we can define an operator
\begin{equation}
\opa_b~:~\cc_{0,q}^s(\mb)\to \cc_{0,q+1}^{s-1}(\mb),
\end{equation}
called the \emph{tangential Cauchy-Riemann operator}, by setting
$\opa_b f=df_{|_{H_{0,1}\mb\times\dots\times H_{0,1}\mb}}$. It satisfies
$\opa_b\circ\opa_b=0$.

\begin{defin}\label{plongement}
Let $(\mb, H_{0,1}\mb)$ be an abstract CR manifold, $X$ be  a
complex manifold  and $F~:~\mb\to X$ be an embedding of class
$\cc^l$, then $F$ is called a \emph{CR embedding} if
$dF(H_{0,1}\mb)$ is a subbundle of the bundle $T_{0,1}X$ of the
antiholomorphic vector fields of $X$ and
$dF(H_{0,1}\mb)=T_{0,1}X\cap T_\cb F(\mb)$.
\end{defin}

Let $F$ be a CR embedding of an abstract CR manifold $\mb$ into a complex
manifold $X$ and set $M=F(\mb)$, then $M$ is a CR manifold with the
CR structure $H_{0,1}M=T_{0,1}X\cap T_\cb M$.

Let $U$ be a coordinate domain in $X$, then
$F_{|_{F^{-1}(U)}}=(f_1,\dots ,f_N)$, with $N=\rm{dim}_\cb X$, and
$F$ is a CR embedding if and only if, for all $1\leq j\leq N$,
$\opa_b f_j=0$.

A CR embedding is called \emph{generic} if $\rm{dim}_\cb~X-
\rm{rk}_\cb~H_{0,1}M=\rm{codim}_\rb~M$.

\begin{defin}\label{perturbation}
An almost CR structure $\widehat H_{0,1}\mb$ on $\mb$ is said to be of
\emph{finite distance} to a given CR structure $H_{0,1}\mb$ if $\widehat
H_{0,1}\mb$ can be represented as a graph in $T_\cb\mb$ over
$H_{0,1}\mb$.

It is called an \emph{horizontal perturbation} of the CR structure
$H_{0,1}\mb$ if it is  of finite distance to  $H_{0,1}\mb$  and
moreover there exists $\Phi\in\cc_{0,1}(\mb,H_{1,0}\mb)$ such that
\begin{equation}\label{perturbation horizontale}
\widehat H_{0,1}\mb=\{\ol W\in T_\cb\mb~|~\ol W=\ol Z-\Phi(\ol Z), \ol
Z\in H_{0,1}\mb\},
\end{equation}
which means that $\widehat H_{0,1}\mb$ is a graph in $H\mb=H_{1,0}\mb\oplus
H_{0,1}\mb$ over $H_{0,1}\mb$.
\end{defin}

Note that an horizontal perturbation of the original CR structure
preserves the complex tangent bundle $H\mb$.

Assume $\mb$ is an abstract CR manifold  and $\widehat H_{0,1}\mb$
is an integrable horizontal perturbation of the original CR
structure $H_{0,1}\mb$ on $\mb$. If $\opa_b^\Phi$ denotes the
tangential Cauchy-Riemann operator associated to the CR structure
$\widehat H_{0,1}\mb$, then we have
\begin{equation}
\opa_b^\Phi=\opa_b-\Phi\lrcorner d=\opa_b-\Phi\lrcorner\pa_b,
\end{equation}
where $\opa_b$ is the tangential Cauchy-Riemann operator associated
to the original CR structure $H_{0,1}\mb$ and $\pa_b$ involves only
holomorphic tangent vector fields.

The annihilator $H^0\mb$ of $H\mb=H_{1,0}\mb\oplus H_{0,1}\mb$ in
$T_\cb^*\mb$ is called the \emph{characteristic bundle} of $\mb$.
Given $p\in\mb$, $\omega\in H^0_p\mb$ and $X,Y\in H_p\mb$, we choose
$\wt\omega\in\Gamma(\mb,H^0\mb)$ and $\wt X,\wt Y\in
\Gamma(\mb,H\mb)$ with $\wt\omega_p=\omega$, $\wt X_p=X$ and $\wt
Y_p=Y$. Then $d\wt\omega(X,Y)=-\omega([\wt X,\wt Y])$. Therefore we
can associate to each $\omega\in H^0_p\mb$ an hermitian form
\begin{equation}\label{levi}
L_\omega(X)=-i\omega([\wt X,\ol{\wt X}])
\end{equation}
on $H_p\mb$. This is called the \emph{Levi form} of $\mb$ at
$\omega\in H^0_p\mb$.

In the study of the $\opa_b$-complex two important geometric
conditions were introduced for CR manifolds of real dimension $2n-k$
and CR-dimension $n-k$. The first one by Kohn in the hypersurface
case, $k=1$, the condition Y(q), the second one by Henkin in
codimension $k$, $k\ge 1$, the $q$-concavity.

An abstract CR manifold $\mb$ of hypersurface type satisfies Kohn's condition $Y(q)$ at a point
$p\in\mb$ for some $0\le q\le n-1$, if the Levi form of $\mb$ at $p$
has at least $\max(n-q,q+1)$ eigenvalues of the same sign or at
least $\min(n-q,q+1)$ eigenvalues of opposite signs.

An abstract CR manifold $\mb$ is said to be \emph{$q$-concave} at $p\in\mb$
for some $0\le q\le n-k$, if the Levi form $L_\omega$ at $\omega\in
H^0_p\mb$ has at least $q$ negative eigenvalues on $H_p\mb$ for
every nonzero $\omega\in H^0_p\mb$.

In \cite{ShWa} the condition Y(q) is extended to arbitrary
codimension.
\begin{defin}\label{y(q)}
An abstract CR manifold is said to satisfy \emph{condition Y(q)} for
some $1\le q\le n-k$ at $p\in\mb$ if the Levi form $L_\omega$ at
$\omega\in H^0_p\mb$ has at least $n-k-q+1$ positive eigenvalues or
at least $q+1$ negative eigenvalues on $H_p\mb$ for every nonzero
$\omega\in H^0_p\mb$.
\end{defin}

Note that in the hypersurface case, i.e. $k=1$, this condition is
equivalent to the classical condition Y(q) of Kohn for
hypersurfaces. Moreover, if $\mb$ is $q$-concave at $p\in\mb$, then
$q\le (n-k)/2$ and condition Y(r) is satisfied at $p\in\mb$ for any
$0\le r\le q-1$ and $n-k-q+1\le r\le n-k$.

%% file: anis3.tex
\section{Stability of  vanishing theorems by horizontal perturbations of the CR structure}\label{s3}

Let $(\mb,H_{0,1}\mb)$ be an abstract compact CR manifold of class
$\ci$, of real dimension $2n-k$ and CR dimension $n-k$, and
$\widehat H_{0,1}\mb$ be an integrable horizontal perturbation of
$H_{0,1}\mb$. We denote by $\mb$ the abstract CR manifold
$(\mb,H_{0,1}\mb)$ and by $\widehat\mb$ the abstract CR manifold
$(\mb,\widehat H_{0,1}\mb)$.

Since $\widehat H_{0,1}\mb$ is an horizontal perturbation of
$H_{0,1}\mb$, which means that $\widehat H_{0,1}\mb$ is a graph in
$H\mb=H_{1,0}\mb\oplus H_{0,1}\mb$ over $H_{0,1}\mb$, the space
$H\widehat\mb=H_{1,0}\widehat\mb\oplus H_{0,1}\widehat\mb$ coincides
with the space $H\mb$ and consequently the two abstract CR manifolds
$\mb$ and $\widehat\mb$ have the same characteristic bundle and
hence the same Levi form. This implies in particular that if $\mb$
satisfies condition Y(q) at each point, then $\widehat\mb$ satisfies
also condition Y(q) at each point.

It follows from the Hodge decomposition theorem and the results in
\cite{ShWa} that if $\mb$ is an abstract compact CR manifold of
class $\ci$ which satisfies condition Y(q) at each point, then the
cohomology groups $H^{p,q}(\mb)$, $0\leq p\leq n$, are finite
dimensional. A natural question is then the stability by small
horizontal perturbations of the CR structure of the vanishing of
these groups.

Let us consider a sequence $(\bc^l(\mb),~l\in\nb)$ of Banach
spaces with
$\bc^{l+1}(\mb)\subset\bc^l(\mb)$, which are invariant by horizontal
perturbations of the CR structure of $\mb$ and such that if
$f\in\bc^l(\mb)$, $l\ge 1$, then $X_\cb f\in\bc^{l-1}(\mb)$ for all
complex tangent vector fields $X_\cb$ to $\mb$ and there exists
$\theta(l)\in\nb$ with $\theta(l+1)\geq\theta(l)$ such that $fg\in\bc^l(\mb)$ 
if $f\in\bc^l(\mb)$ and 
$g\in\cc^{\theta(l)}(\mb)$. Such a sequence $(\bc^l(\mb),~l\in\nb)$
will be called \emph{a sequence of anisotropic spaces}. We denote by
$\bc^l_{p,r}(\mb)$ the space of $(p,r)$-forms on $\mb$ whose
coefficients belong to $\bc^l(\mb)$. Moreover we will say that these
Banach spaces are \emph{adapted to the $\opa_b$-equation in degree
$r\geq 1$} if, when $H^{p,r}(\mb)=0$, $0\le p\le n$, there exist
linear continuous operators $A_s$, $s=r,r+1$, from $\bc^0_{p,s}(\mb)$
into $\bc^0_{p,s-1}(\mb)$ which are also continuous from
$\bc^l_{p,s}(\mb)$ into $\bc^{l+1}_{p,s-1}(\mb)$, $l\in\nb$, and
moreover satisfy
\begin{equation}\label{homotopieglobalebis}
f=\opa_b A_r f+A_{r+1}\opa_b f\,,
\end{equation}
for $f\in\bc^l_{p,r}(\mb)$ with $\opa_b f=0$.

\begin{thm}\label{annulation}
Let $\mb=(\mb,H_{0,1}\mb)$ be an abstract compact CR manifold of
class $\ci$, of real dimension $2n-k$ and CR dimension $n-k$, and
$\widehat\mb=(\mb,\widehat H_{0,1}\mb)$ another abstract compact CR
manifold such that $\widehat H_{0,1}\mb$ is an integrable horizontal
perturbation of $H_{0,1}\mb$. Let also $(\bc^l(\mb),~l\in\nb)$ be a
sequence of anisotropic Banach spaces and $q$ be an integer, $1\le
q\le (n-k)/2$. Finally let
$\Phi\in\cc^{\theta(l)}_{0,1}(\mb,H_{1,0}\mb)$ be the differential
form which defines the tangential Cauchy-Riemann operator
$\opa_b^\Phi=\opa_b-\Phi\lrcorner\pa_b$ associated to the CR
structure $\widehat H_{0,1}\mb$.

Assume $H^{p,r-p}(\mb)=0$, for $1\leq p\leq r$ and $1\leq r\leq
s_1(q)$, or $H^{n-p,r+p}(\mb)=0$, for $0\le p\le n-k-r$ and
$s_2(q)\le r\le n-k$ and that the Banach spaces $(\bc^l(\mb),~l\in\nb)$
are adapted to the $\opa_b$-equation in degree $r$, $1\leq r\leq
s_1(q)$ or $s_2(q)\le r\le n-k$. Then, for each $l\in\nb$, there
exists $\delta>0$ such that, if $\|\Phi\|_{\theta(l)}<\delta$,

(i) for each $\opa_b^\phi$-closed form $f$ in
$\bc^l_{0,r}(\widehat\mb)$, $1\le r\le s_1(q)$, such that the part
of $f$ of bidegree $(0,r)$ for the initial CR structure $H_{0,1}\mb$
vanishes, there exists a form $u$ in
$\bc^{l+1}_{0,r-1}(\widehat\mb)$ satisfying $\opa_b^\phi u=f$,

(ii) for each $\opa_b^\phi$-closed form $f$ in
$\bc^l_{n,r}(\widehat\mb)$, $s_2(q)\le r\le n-k$, there exists a
form $u$ in $\bc^{l+1}_{n,r-1}(\widehat\mb)$ satisfying $\opa_b^\phi
u=f$.
\end{thm}

\begin{rem}
Note that if both $\mb=(\mb,H_{0,1}\mb)$ and
$\widehat\mb=(\mb,\widehat H_{0,1}\mb)$ are embeddable in the same
complex manifold $X$, any $r$-form on the differential manifold
$\mb$, which represents a form of bidegree $(0,r)$ for the CR
structure $\widehat H_{0,1}\mb$ represents also a form of bidegree
$(0,r)$ for the CR structure $H_{0,1}\mb$. Hence the bidegree
hypothesis in (i) of Theorem \ref{annulation} is automatically
fulfilled.
\end{rem}

\begin{proof}
Let $f\in\bc^l_{0,r}(\widehat\mb)$ be a $(0,r)$-form for the CR
structure $\widehat H_{0,1}\mb$, $1\le r\le s_1(q)$, such that
$\opa_b^\Phi f=0$, we want to solve the equation
\begin{equation}\label{dbarphi0}
\opa_b^\Phi u=f.
\end{equation}

The form $f$ can be written $\sum_{p=0}^r f_{p,r-p}$, where the
forms $f_{p,r-p}$ are of type $(p,r-p)$ for the CR structure
$H_{0,1}\mb$. Then by considerations of bidegrees, the equation
$\opa_b^\Phi f=0$ is equivalent to the family of equations
$\opa_b^\Phi f_{p,r-p}=0$, $0\le p\le r$.

Moreover, if $u=\sum_{s=0}^{r-1} u_{s,r-1-s}$, where the forms
$u_{s,r-1-s}$ are of type $(s,r-1-s)$ for the CR structure
$H_{0,1}\mb$, is a solution of \eqref{dbarphi0}, then
\begin{equation*}
\opa_b^\Phi u_{p,r-1-p}=f_{p,r-p}\,,
\end{equation*}
for $0\le p\le r-1$, and
$$f_{r,0}=0\,.$$

Therefore a necessary condition on $f$ for the solvability of
\eqref{dbarphi0} is that $\opa_b^\Phi f=0$ and $f_{r,0}=0$, where
$f_{r,0}$ is the part of type $(r,0)$ of $f$ for the CR structure
$H_{0,1}\mb$, and, to solve \eqref{dbarphi0}, we have to consider
the equation
\begin{equation}\label{dbarphi}
\opa_b^\Phi v=g,
\end{equation}
where $g\in\bc^l_{p,r-p}(\mb)$ is a $(p,r-p)$-form for the CR
structure $H_{0,1}\mb$, $0\le p\le r-1$, which is $\opa_b^\Phi$
closed. By definition of the operator $\opa_b^\Phi$, this means
solving the equation $\opa_b v=g+\Phi\lrcorner\pa_b v$. Consequently
if $v$ is a solution of (\ref{dbarphi}), then $\opa_b
(g+\Phi\lrcorner\pa_b v)=0$ and by (\ref{homotopieglobalebis})
$$\opa_b(A_{r-p}(g+\Phi\lrcorner\pa_b v))=g+\Phi\lrcorner\pa_b v.$$

Assume $\Phi$ is of class $\cc^{\theta(l)}$, then the map
\begin{align*}
\Theta~:~\bc^{l+1}_{p,r-1}(\mb)&\to\bc^{l+1}_{p,r-1}(\mb)
\\
v&\mapsto A_{r-p} g+A_{r-p}(\Phi\lrcorner\pa_b v)\,.
\end{align*}
is continuous, and the fixed points of $\Theta$ are good candidates
to be solutions of (\ref{dbarphi}).

Let $\delta_0$ such that, if $\|\Phi\|_{\theta(l)}<\delta_0$, then
the norm of the bounded endomorphism
$A_{r-p}\circ\Phi\lrcorner\pa_b$ of $\bc^{l+1}_{p,r-p-1}(\mb)$ is
equal to $\epsilon_0<1$. We shall prove that, if
$\|\Phi\|_{\theta(l)}<\delta_0$, $\Theta$ admits a unique fixed
point.

Consider first the uniqueness of the fixed point. Assume $v_1$ and
$v_2$ are two fixed points of $\Theta$, then
\begin{align*}
v_1&=\Theta(v_1)=A_{r-p} g+A_{r-p}(\Phi\lrcorner\pa_b v_1)\\
v_2&=\Theta(v_2)=A_{r-p} g+A_{r-p}(\Phi\lrcorner\pa_b v_2).
\end{align*}
This implies
$$v_1-v_2=A_{r-p}\big(\Phi\lrcorner\pa_b (v_1-v_2)\big)$$
and, by the hypothesis on $\Phi$,
$$\|v_1-v_2\|_{\bc^{l+1}}<\|v_1-v_2\|_{\bc^{l+1}}\quad{\rm or}\quad v_1=v_2$$
 and hence $v_1=v_2$.

For the existence we proceed by iteration. We set
$v_0=\Theta(0)=A_{r-p}(g)$ and, for $n\geq 0$,
$v_{n+1}=\Theta(v_n)$. Then for $n\geq 0$, we get
$$v_{n+1}-v_n=A_{r-p}(\Phi\lrcorner\pa_b (v_n-v_{n-1})).$$

Therefore, if $\|\Phi\|_{\theta(l)}<\delta_0$, the sequence
$(v_n)_{n\in\nb}$ is a Cauchy sequence in the Banach space
$\bc^{l+1}_{p,r-1}(\mb)$ and hence converges to a form $v$,
moreover by continuity of the map $\Theta$, $v$ satisfies
$\Theta(v)=v$.

It remains to prove that $v$ is a solution of (\ref{dbarphi}). Since
$H^{p,r-p}(\mb)=0$ for $1\le p\le r$, it follows from
(\ref{homotopieglobalebis}) and from the definition of the sequence
$(v_n)_{n\in\nb}$ that
\begin{equation*}
g-\opa_b^\Phi
v_{n+1}=\Phi\lrcorner\pa_b(v_{n+1}-v_n)+A_{r-p+1}\opa_b(g+\Phi\lrcorner\pa_b
v_n)
\end{equation*}
and since
\begin{align*}
\opa_b(g+\Phi\lrcorner\pa_b v_n)&=\opa_b g-\opa_b(\opa_b-\Phi\lrcorner\pa_b) v_n\\
&=\opa_b g-\opa_b(\opa_b^\Phi v_n)\\
&=\opa_b g-(\opa_b^\Phi+\Phi\lrcorner\pa_b)(\opa_b^\Phi v_n)\\
&=\opa_b g-\Phi\lrcorner\pa_b(\opa_b^\Phi v_n), \quad
\text{since}~(\opa_b^\Phi)^2=0\\
&=\Phi\lrcorner\pa_b(g-\opa_b^\Phi v_n), \quad
\text{since}~\opa_b^\Phi g=0,
\end{align*}
we get
\begin{equation}\label{erreur}
g-\opa_b^\Phi
v_{n+1}=\Phi\lrcorner\pa_b(v_{n+1}-v_n)+A_{r-p+1}(\Phi\lrcorner\pa_b(g-\opa_b^\Phi
v_n)).
\end{equation}

Note that since $g\in\bc^l_{p,r}(\mb)$ and $\Phi$ is of
class $\cc^{\theta(l)}$, it follows from the definition of the $v_n$s 
that $v_n\in\bc^{l+1}_{p,r}(\mb)$ and $\opa_b^\Phi
v_n\in\bc^l_{p,r}(\mb)$ for all $n\in\nb$.

Thus by (\ref{erreur}), we have the
estimate
\begin{equation}\label{est1}
\|g-\opa_b^\Phi v_{n+1}\|_{\bc^l}
\leq\|\Phi\lrcorner\pa_b\|\|(v_{n+1}-v_n)\|_{\bc^{l+1}}+\|A_{r-p+1}\circ\Phi\lrcorner\pa_b\|\|g-\opa_b^\Phi
v_n\|_{\bc^l}.
\end{equation}
Let $\delta$ such that if $\|\Phi\|_{\theta(l)}<\delta$, then the
maximum of the norm of the bounded endomorphisms
$A_s\circ\Phi\lrcorner\pa_b$, $s=r-p,r-p+1$, of
$\bc^{l}_{p,s-1}(\mb)$ is equal to $\epsilon<1$. Assume
$\|\Phi\|_{\theta(l)}<\delta$, then by induction we get
\begin{equation}\label{est2}
\|g-\opa_b^\Phi v_{n+1}\|_{\bc^l} \leq
(n+1)\epsilon^{n+1}\|\Phi\lrcorner\pa_b\|\|v_0\|_{\bc^{l+1}}+\epsilon^{n+1}\|g-\opa_b^\Phi
v_0\|_{\bc^l}.
\end{equation}
But $g-\opa_b^\Phi v_0=\Phi\lrcorner\pa_b
A_{r-p}g+A_{r-p+1}(\Phi\lrcorner\pa_b g)$ and hence $\|g-\opa_b^\Phi
v_0\|_{\bc^l}\leq\|\Phi\lrcorner\pa_b\|\|A_{r-p}g\|_{\bc^{l+1}}
+\epsilon \|g\|_{\bc^l}$. This implies
\begin{equation}\label{est3}
\|g-\opa_b^\Phi v_{n+1}\|_{\bc^l}\leq
(n+2)\epsilon^{n+1}\|\Phi\lrcorner\pa_b\|\|A_{r-p}\|\|g\|_{\bc^{l}}+\epsilon^{n+2}\|g\|_{\bc^l}.
\end{equation}
Since $\epsilon<1$, the righthand side of (\ref{est3}) tends to
zero, when $n$ tends to infinity and by continuity of the operator
$\opa_b^\Phi$ from $\bc^{l+1}_{p,r-1}(\mb)$ into
$\bc^l_{p,r}(\mb)$, the lefthand side of (\ref{est3}) tends
to $\|g-\opa_b^\Phi v\|_{\bc^l}$, when $n$ tends to infinity, which
implies that $v$ is a solution of (\ref{dbarphi}).

Now if $f\in\bc^l_{n,r}(\widehat\mb)$ is an $(n,r)$-form for the CR
structure $\widehat H_{0,1}\mb$, $s_2(q)\le r\le n-k$, such that
$\opa_b^\Phi f=0$, then the form $f$ can be written
$\sum_{p=0}^{n-k-r} f_{n-p,r+p}$, where the forms $f_{n-p,r+p}$ are
$\opa_b^\Phi$-closed and of type $(n-p,r+p)$ for the CR structure
$H_{0,1}\mb$. Then to solve the equation
\begin{equation*}
\opa_b^\Phi u=f,
\end{equation*}
it is sufficient to solve the equation $\opa_b^\Phi v=g$ for
$g\in\bc^l_{n-p,r+p}(\mb)$ and this can be done in the same way as
in the case of the small degrees, but using the vanishing of the
cohomology groups $H^{n-p,p+r}(\mb)$ for $0\le p\le n-k-r$ and
$s_2(q)\le r\le n-k$.
\end{proof}

Assume the horizontal perturbation of the original CR structure on
$\mb$ is smooth, i.e. $\Phi$ is of class $\ci$, then we can  defined
on $\widehat\mb$ the cohomology groups
$$H^{o,r}_\Phi(\widehat\mb)=\{f\in\ci_{o,r}(\widehat\mb)~|~f_{r,0}=0,~\opa_b^\Phi f=0\}/\opa_b^\Phi (\ci_{o,r-1}(\widehat\mb))$$
and
$$H^{n,r}_\Phi(\widehat\mb)=\{f\in\ci_{o,r}(\widehat\mb)~|~\opa_b^\Phi f=0\}/ \opa_b^\Phi (\ci_{n,r-1}(\widehat\mb))$$
for $1\le r\le n-k$.

\begin{cor}\label{lisse}
Under the hypotheses of Theorem \ref{annulation}, if the sequence
$(\bc^l(\mb),~l\in\nb)$ of anisotropic Banach spaces is such that
$\cap_{l\in\nb} B_l(\mb)=\ci(\mb)$ and if the horizontal
perturbation of the original CR structure on $\mb$ is smooth, there
exists a sequence $(\delta_l)_{l\in\nb}$ of positive real numbers
such that, if $\|\Phi\|_{\theta(l)}<\delta_l$ for each $l\in\nb$

(i)  $H^{p,r-p}(\mb)=0$, for all $1\leq p\leq r$, implies
$H^{o,r}_\Phi(\widehat\mb)=0$, when $1\leq r\leq s_1(q)$,

(ii) $H^{n-p,r+p}(\mb)=0$, for all $0\le p\le n-k-r$, implies
$H^{n,r}_\Phi(\widehat\mb)=0$, when $s_2(q)\le r\le n-k$.
\end{cor}
\begin{proof}
It is a direct consequence of the proof of Theorem \ref{annulation}
by the uniqueness of the fixed point of $\Theta$.
\end{proof}

%%% Local Variables:
%%% mode: latex
%%% TeX-master: "anisotropic"
%%% End:

%% file: anis4.tex
\section{Anisotropic spaces}\label{s4}

In the previous section the main theorem is proved under the
assumption of the existence of sequences of anisotropic spaces on
abstract CR manifolds satisfying good properties with respect to the
tangential Cauchy-Riemann operator. We will precise Theorem
\ref{annulation} by considering some Sobolev and some H\"older
anisotropic spaces for which global homotopy formulas for the
tangential Cauchy-Riemann equation with good estimates hold
under some geometrical conditions.

In this section $\mb=(\mb,H_{0,1}\mb)$ denotes an abstract compact
CR manifold of class $\ci$, of real dimension $2n-k$ and CR
dimension $n-k$.

Let us define some anisotropic Sobolev spaces of functions:

- $\cs^{0,p}(\mb)$, $1<p<\infty$, is the set of $L^p_{loc}$
functions on $\mb$.

- $\cs^{1,p}(\mb)$, $1<p<\infty$, is the set of fonctions on $\mb$
such that $f\in W^{\frac{1}{2},p}(\mb)$ and $X_\cb f\in
L^p_{loc}(\mb)$, for all complex tangent vector fields $X_\cb$ to
$M$.

- $\cs^{l,p}(\mb)$, $l\geq 2$, $1<p<\infty$, is the set of functions
$f$ such that $Xf\in\cs^{l-2,p}(\mb)$, for all tangent vector fields
$X$ to $M$ and $X_\cb f\in\ac^{l-1,p}(\mb)$, for all complex tangent
vector fields $X_\cb$ to $M$.

The sequence $(\cs^{l,p}(\mb),~l\in\nb)$ is a sequence of
anisotropic spaces in the sense of section \ref{s3} with
$\theta(l)=l+1$. Moreover $\cap_{l\in\nb}\cs^{l,p}(\mb)=\ci(\mb)$.

The anisotropic H\"older space of forms $\cs^{l,p}_*(\mb)$, $l\geq
0$, $1<p<\infty$, is then the space of forms on $\mb$, whose
coefficients are in  $\cs^{l,p}(\mb)$.

We have now to see if the sequence $(\cs^{l,p}(\mb),~l\in\nb)$ is
adapted to the $\opa_b$-equation for some degree $r$.

The $L^2$ theory for $\Box_b$ in abstract CR manifolds of arbitrary
codimension is developed in \cite{ShWa}. There it is proved that if $\mb$
satisfies condition Y(r) the Hodge decomposition theorem holds in
degree $r$, which means that there exist a compact operator
$N_b~:~L^2_{p,r}(\mb)\to Dom(\Box_b)$ and a continuous operator
$H_b~:~L^2_{p,r}(\mb)\to L^2_{p,r}(\mb)$ such that for any $f\in
L^2_{p,r}(\mb)$
\begin{equation}\label{hodge}
f=\opa_b\opa_b^* N_b f+\opa_b^*\opa_b N_b f+H_b f
\end{equation}
Moreover $H_b$ vanishes on exact forms and if $N_b$ is also defined
on $L^2_{p,r+1}(\mb)$ then $N_b\opa_b=\opa_b N_b$.

Therefore if $\mb$ satisfy both conditions Y(r) and Y(r+1) then
\eqref{hodge} becomes an homotopy formula and using the Sobolev and
the anisotropic Sobolev estimates in \cite{ShWa} (Theorems 3.3 and
Corollary 1.3 (2)) we get the following result:

\begin{prop}\label{sob}
If $\mb$ is $q$-concave, the sequence $(\cs^{l,p}(\mb),~l\in\nb)$ of
anisotropic spaces is adapted to the $\opa_b$-equation in degree $r$
for $0\leq r\leq q-2$ and $n-k-q+1\le r\le n-k$.
\end{prop}

Let us define now some anisotropic H\"older spaces of functions:

- $\ac^\alpha(\mb)$, $0<\alpha<1$, is the set of continuous
functions on $\mb$ which are in $\cc^{\alpha/2}(\mb)$.

- $\ac^{1+\alpha}(\mb)$, $0<\alpha<1$, is the set of functions $f$
such that $f\in\cc^{(1+\alpha)/2}(\mb)$ and $X_\cb
f\in\cc^{\alpha/2}(\mb)$, for all complex tangent vector fields
$X_\cb$ to $\mb$.  Set
\begin{equation}\label{holderanis}
  \|f\|_{A\alpha}=\|f\|_{(1+\alpha)/2}+\sup_{\|X_\cb\|\leq 1}\|X_\cb f\|_{\alpha/2}
\end{equation}

- $\ac^{l+\alpha}(\mb)$, $l\geq 2$, $0<\alpha<1$, is the set of
functions $f$ of class $\cc^{[l/2]}$ such that
$Xf\in\ac^{l-2+\alpha}(\mb)$, for all tangent vector fields $X$ to
$M$ and $X_\cb f\in\ac^{l-1+\alpha}(\mb)$, for all complex tangent
vector fields $X_\cb$ to $\mb$.

The sequence $(\ac^{l+\alpha}(\mb),~l\in\nb)$ is a sequence of
anisotropic spaces in the sense of section \ref{s3} with
$\theta(l)=l+1$. Moreover
$\cap_{l\in\nb}\ac^{l+\alpha}(\mb)=\ci(\mb)$.

The anisotropic H\"older space of forms $\ac^{l+\alpha}_*(\mb)$,
$l\geq 0$, $0<\alpha<1$, is then the space of continuous forms on
$\mb$, whose coefficients are in  $\ac^{l+\alpha}(\mb)$.

It remains to see if the sequence $(\ac^{l+\alpha}(\mb),~l\in\nb)$
is adapted to the $\opa_b$-equation for some degrees $r$.

Assume $\mb$ is locally embeddable and $1$-concave. Then, by
Proposition 3.1 in \cite{HiNa1}, there exist a complex manifold $X$
and a smooth generic embedding $\ec~:~\mb\to M\subset X$ such that
$M$ is a smooth compact CR submanifold of $X$ with the CR structure
$H_{0,1}M=d\ec(H_{0,1}\mb)=T_\cb M\cap T_{0,1}$. If $E$ is a CR
vector bundle over $\mb$, by the $1$-concavity of $\mb$ and after an
identification between $\mb$ and $M$, the CR bundle $E$ can be
extended to an holomorphic bundle in a neighborhood of $M$, which we
still denote by $E$. With these notations it follows from
\cite{Laanis} that if $\mb$ is $q$-concave, $q\ge 1$, there exist
finite dimensional subspaces $\Cal H_r$ of $\Cal
Z^\infty_{n,r}(\mb,E)$, $0\le r\le q-1$ and $n-k-q+1\le r\le n-k$,
where $\Cal H_0=\Cal Z^\infty_{n,0}(\mb,E)$, continuous linear
operators
\begin{equation*}
A_r~:~\cc^0_{n,r}(\mb,E)\to \cc^0_{n,r-1}(\mb,E),\quad 1\leq r\leq
q~ {\rm and}~ n-k-q+1\le r\le n-k
\end{equation*}
and continuous linear projections $$ P_r:\Cal
C^0_{n,r}(\mb,E)\rightarrow \Cal C^0_{n,r}(\mb,E) \,,\qquad 0\le
r\le q-1~ {\rm and}~ n-k-q+1\le r\le n-k\,,
$$with
\begin{equation}\label{19.4.07'}
\im P_r=\Cal H_r\,,\qquad 0\le r\le q-1~ {\rm and}~ n-k-q+1\le r\le
n-k\,,
\end{equation}and
\begin{equation}\label{20.4.07}
\cc^0_{n,r}(M,E)\cap \opa_b\cc^0_{n,r-1}(M,E)\subseteq\ke P_r\,,\qquad 1\le
r\le q-1~ {\rm and}~ n-k-q+1\le r\le n-k,
\end{equation} such that:

(i) For all $l\in\nb$ and $1\leq r\leq q$ or $n-k-q+1\le r\le n-k$,
\begin{equation*}
A_r(\ac^{l+\alpha}_{n,r}(\mb,E))\subset
\ac^{l+1¨\alpha}_{n,r-1}(\mb,E)
\end{equation*}
and $A_r$ is continuous as an operator between
$\ac^{l+\alpha}_{n,r}(\mb,E)$ and $\ac^{l+1+\alpha}_{n,r-1}(\mb,E)$.

(ii) For all $0\leq r\leq q-1$ or $n-k-q+1\le r\le n-k$ and
$f\in\cc^0_{n,r}(\mb,E)$ with $\opa_bf\in\cc^0_{n,r+1}(\mb,E)$,
\begin{equation}\label{homotopieglobale}
f-P_rf=\begin{cases}A_{1}\opa_b f&\text{if }r=0\,,\\\opa_b
A_rf+A_{r+1}\opa_b f\qquad&\text{if }1\leq r\leq q-1 ~ {\rm or}~
n-k-q+1\le r\le n-k \,.
\end{cases}
\end{equation}

This implies the following result
\begin{prop}\label{hold}
If $\mb$ is locally embeddable and $q$-concave the sequence
$\ac^{l+\alpha}(\mb),~l\in\nb$ of anisotropic spaces is adapted to
the $\opa_b$-equation in degree $r$ for $1\leq r\leq q-1$ and
$n-k-q+1\le r\le n-k$
\end{prop}

Finally let us recall the definition of the anisotropic H\"older
spaces $\Gamma^{l+\alpha}(\mb)$ of Folland and Stein.

- $\Gamma^\alpha(\mb)$, $0<\alpha<1$, is the set of continuous
fonctions in $\mb$ such that if for every $x_0\in\mb$
$$\sup_{\gamma(.)}\frac{|f(\gamma(t)-f(x_0)|}{|t|^\alpha}<\infty$$
for any complex tangent curve $\gamma$ through $x_0$.

- $\Gamma^{l+\alpha}(\mb)$, $l\ge 1$, $0<\alpha<1$, is the set of
continuous fonctions in $\mb$ such that $X_\cb
f\in\Gamma^{l-1+\alpha}(\mb)$, for all complex tangent vector fields
$X_\cb$ to $\mb$.

The spaces $\Gamma^{l+\alpha}(\mb)$ are subspaces of the spaces 
$\ac^{l+\alpha}(\mb)$.

Note that by Corollary 1.3 (1) in \cite{ShWa} and Section 3 in
\cite{Laanis} Propositions \ref{sob} and \ref{hold} hold also for
the anisotropic H\"older spaces $\Gamma^{l+\alpha}(\mb)$ of Folland
and Stein.

Let us summarize all this in connection with section \ref{s3} in the
next theorem.

\begin{thm}\label{annul}
If $\mb$ is $q$-concave,

(i) Theorem \ref{annulation} holds  for $\bc^l(\mb)=\cs^{l,p}(\mb)$
with $s_1(q)=q-2$ and $s_2(q)=n-k-q+1$ in the abstract case,

(ii) Theorem \ref{annulation} holds  for
$\bc^l(\mb)=\ac^{l+\alpha}(\mb)$ with $s_1(q)=q-1$ and
$s_2(q)=n-k-q+1$ when $\mb$ is locally embeddable

(iii) Theorem \ref{annulation} holds  for
$\bc^l(\mb)=\Gamma^{l+\alpha}(\mb)$ with $s_1(q)=q-2$ in the
abstract case and $s_1(q)=q-1$ when $\mb$ is locally embeddable, and
with $s_2(q)=n-k-q+1$ in both cases.
\end{thm}

Since in all the three cases of Theorem \ref{annul} we have
$\cap_{l\in\nb} B_l(\mb)=\ci(\mb)$, Corollary \ref{lisse} becomes

\begin{cor}
Let $\mb=(\mb,H_{0,1}\mb)$ be an abstract compact CR manifold of
class $\ci$, of real dimension $2n-k$ and CR dimension $n-k$, and
$\widehat\mb=(\mb,\widehat H_{0,1}\mb)$ another abstract compact CR
manifold such that $\widehat H_{0,1}\mb$ is an integrable horizontal
smooth perturbation of $H_{0,1}\mb$.  Let
$\Phi\in\ci_{0,1}(\mb,H_{1,0}\mb)$ be the differential form which
defines the tangential Cauchy-Riemann operator
$\opa_b^\Phi=\opa_b-\Phi\lrcorner\pa_b$ associated to the CR
structure $\widehat H_{0,1}\mb$. Assume $\mb$ is $q$-concave, then
there exists a sequence $(\delta_l)_{l\in\nb}$ of positive real
numbers such that, if $\|\Phi\|_{l}<\delta_l$ for each $l\in\nb$

(i)  $H^{p,r-p}(\mb)=0$, for all $1\leq p\leq r$, implies
$H^{0,r}_\Phi(\widehat\mb)=0$, when $1\leq r\leq q-2$ in the
abstract case and also for $r=q-1$ if $\mb$ is locally embeddable,

(ii) $H^{n-p,r+p}(\mb)=0$, for all $0\le p\le n-k-r$, implies
$H^{n,r}_\Phi(\widehat\mb)=0$, when $n-k-q+1\le r\le n-k$.
\end{cor}

%% file: stability.bbl
\providecommand{\bysame}{\leavevmode\hbox to3em{\hrulefill}\thinspace}
\providecommand{\MR}{\relax\ifhmode\unskip\space\fi MR }
\providecommand{\MRhref}[2]{%
  \href{http://www.ams.org/mathscinet-getitem?mr=#1}{#2}
}
\providecommand{\href}[2]{#2}
\begin{thebibliography}{1}

\bibitem{BaLa}
M.~Y. Barkatou and C.~Laurent-Thi{\'e}baut, \emph{Estimations optimales pour
  l'op{\'e}rateur de {C}auchy-{R}iemann tangentiel}, Michigan Math. Journal
  \textbf{54} (2006), 545--586.

\bibitem{BroLaLe}
T.~Br{\"o}nnle, C.~Laurent-Thi{\'e}baut, and J.~Leiterer, \emph{Global homotopy
  formulas on q-concave {CR} manifolds for large degrees}, Preprint.

\bibitem{FoSt}
G.~B. Folland and E.M. Stein, \emph{Estimates for the $\opa_b$-complex and
  analysis on the heisenberg group}, Comm. Pure Appl. Math. \textbf{27} (1974),
  429--522.

\bibitem{HiNa1}
C.D. Hill and M.~Nacinovich, \emph{Pseudoconcave {C}{R} manifolds}, Complex
  analysis and geometry, V. Ancona, E. Ballico, A. Silva, eds., Marcel Decker,
  Inc., New-York, 1996, pp.~275--297.

\bibitem{Laanis}
C.~Laurent-Thi{\'e}baut, \emph{Poincar{\'e} lemma and global homotopy formulas
  with sharp anisotropic {H}{\"o}lder estimates in q-concave {C}{R} manifolds},
  Pr{\'e}publication de l'Institut Fourier hal-00335229, ar{X}iv 0810.5295
  (2008), 1--15.

\bibitem{LaLeglobal}
C.~Laurent-Thi{\'e}baut and J.~Leiterer, \emph{Global homotopy formulas on
  q-concave {CR} manifolds for small degrees}, J. Geom. Anal. \textbf{18}
  (2008), 511--536.

\bibitem{Poly1}
P.~L. Polyakov, \emph{Global $\opa_{M}$-homotopy with ${\cc}^k$ estimates for a
  family of compact, regular $q$-pseudoconcave {CR} manifolds}, Math.
  Zeitschrift \textbf{247} (2004), 813--862.

\bibitem{Poly2}
\bysame, \emph{Versal embeddings of compact $3$-pseudoconcave
  {CR}-submanifolds}, Math. Zeitschrift \textbf{248} (2004), 267--312.

\bibitem{ShWa}
M.-C. Shaw and L.~Wang, \emph{{H}{\"o}lder and ${L}^p$ estimates for
  $\square_b$ on {CR} manifolds of arbitrary codimension}, Math. Ann.
  \textbf{331} (2005), 297--343.

\end{thebibliography}
